\documentclass[12pt,letterpaper]{amsart}
\usepackage{amssymb,latexsym, amsmath, amsxtra}
\usepackage[dvips]{graphics}

\theoremstyle{plain}
\newtheorem{theorem}{Theorem}[section]

\newtheorem{proposition}[theorem]{Proposition}
\theoremstyle{definition}

\theoremstyle{remark}

\numberwithin{equation}{section}

\newcommand{\Z}{\mathbb Z}
\renewcommand{\H}{\mathbb H}

\def\({\left(}
\def\){\right)}

\def \mm#1#2#3#4 {\begin{pmatrix} \mathstrut#1 & \mathstrut#2 \\
                        \mathstrut#3 & \mathstrut#4\end{pmatrix} }

\begin{document}
\title{Maass forms and their $L$-functions}
\author{David W. Farmer and Stefan Lemurell}

\thanks{Research of the first author supported by the
American Institute of Mathematics
and the NSF Focused Research
Group grant DMS 0244660.
Research of the second author supported in
part by ``Stiftelsen f\"or internationalisering av
h\"ogre utbildning och forskning'' (STINT)
}

\address{
{\parskip 0pt
American Institute of Mathematics\endgraf
360 Portage Ave.\endgraf
Palo Alto, CA 94306\endgraf
farmer@aimath.org\endgraf
\null
Chalmers University of Technology,\endgraf
SE-412 96\endgraf
G\"oteborg, Sweden\endgraf
sj@math.chalmers.se\endgraf
}}

\thispagestyle{empty}
\vspace{.5cm}
\begin{abstract}
We present examples of Maass forms on Hecke congruence groups,
giving low eigenvalues on $\Gamma_0(p)$ for small prime $p$,
and the first 1000 eigenvalues for~$\Gamma_0(11)$.
We also present calculations of the $L$-functions associated to
the Maass forms and make comparisons to the predictions
from random matrix theory.
\end{abstract}

\maketitle

\section{Introduction}

Much recent progress in understanding 
$L$-functions has come from the idea of a ``family'' of
$L$-functions with an associated symmetry type~\cite{KSa}.
The idea is that to a naturally occurring collection of $L$-functions
one can associate a
classical compact group: unitary, symplectic, or orthogonal.
One expects the analytic properties of the $L$-functions
to be largely governed by the symmetry type.
This philosophy has produces a wealth of interesting predictions
which have been confirmed both theoretically and numerically.
See~\cite{CFKRS} for an extensive discussion.

The family of $L$-functions of interest to us here is the
collection of $L$-functions associated to Maass forms.  Specifically,
for a given Hecke congruence group~$\Gamma=\Gamma_0(N)$, 
we consider the Maass forms on $\Gamma$ and the $L$-functions
associated to those Maass forms.  This constitutes an
``orthogonal'' family of $L$-functions, and this leads to 
some specific predictions about statistical properties
of the critical values and the zeros of the $L$-functions.
In this paper we provide the first numerical tests of these
predictions by finding the first 1000 newform Maass forms on
$\Gamma_0(11)$ and computing the associated $L$-functions.
We also test standard predictions about the statistics
of those eigenvalues, and we find the first few eigenvalues
on $\Gamma_0(p)$ for small~$p$.

In the next section we recall properties of Maass forms and their 
associated $L$-functions.  In section~\ref{sec:thealgorithm}
we describe our algorithm for locating and computing Maass forms.
We then present the results of our calculations, in
section~\ref{sec:loweigenvalues} addressing the low eigenvalues
on $\Gamma_0(p)$, and in section~\ref{sec:gamma011} the first
1000 eigenvalues for~$\Gamma_0(11)$.
In section~\ref{sec:examples} we give
example plots of Maass $L$-functions and discuss their general features.
In section~\ref{sec:data} 
we compare the random matrix predictions to our $L$-function data.
Finally, in section~\ref{sec:computing} we describe how the
$L$-functions were computed.

\section{Maass forms and their $L$-Functions}

The $L$-functions we consider are obtained from Maass forms
on the Hecke congruence group $\Gamma_0(11)$.  We first recall the
definition of Maass form, and then describe the connection with 
$L$-functions.  A good reference on Maass forms is
Iwaniec's book~\cite{Iw}.

\subsection{Maass forms}
A {\it Maass form} on
a group $\Gamma\subset PSL(2,\mathbb R)$ is a function $f:{\mathcal H}\to\mathbb R$ 
which satisfies:
\begin{itemize}
\item{}(1.1) $f(\gamma z)=f(z)$ for all $\gamma\in\Gamma$,
\item{}(1.2) $f$ vanishes at the cusps of $\Gamma$, and
\item{}(1.3) $\Delta f = \lambda f$ for some~$\lambda>0$,
\end{itemize}

\noindent where
$$
\Delta=-y^2\left(\frac{\partial^2}{\partial x^2}
        + \frac{\partial^2}{\partial y^2}\right)
$$
is the Laplace--Beltrami operator on~$\mathcal H$.
We set $\lambda=\frac14+R^2$.

In number theory, Maass forms most commonly arise on
Hecke congruence groups:
$$
\Gamma_0(N)=\left\{\mm abcd \in PSL(2,\mathbb Z) \ :\ N|c \right\}.
$$
Here we consider newforms on $\Gamma_0(N)$.  This implies that
\begin{equation}\label{eqn:fourier}
f(z)=\sum_{n=1}^\infty a_n \sqrt y K_{iR} (2\pi n y)SC (2\pi nx),
\end{equation}
where $K_\nu(\cdot)$
is the $K$-Bessel function  
and $SC(x)$  is either $\sin(x)$ or $\cos(x)$. 
In the first case we say
that $f$ is ``odd'' and in the second $f$ is ``even.''
Furthermore, $f$ is an
eigenfunction of the Fricke involution
\begin{equation}\label{eqn:fricke}
f(z)=\pm f\left( -\frac{1}{Nz} \right),
\end{equation}
and $f$ is also a simultaneous eigenfunction of the  Hecke operators
$T_p$ for $p\nmid N$.  A good introduction to this material is \cite{Iw}.
In this paper we only use the properties which are explicitly described above.

If $M|N$ then $\Gamma_0(N)\subset\Gamma_0(M)$.  In particular, if
$f(z)$ is a Maass form on $\Gamma_0(M)$ then $f(kz)$ is a Maass form
on $\Gamma_0(N)$ for all $k|(N/M)$.  Such functions are called
``oldforms'' on $\Gamma_0(N)$ and it is natural to avoid such
functions in a search for Maass forms, for they naturally
belong on the larger group.  The Maass forms which naturally
live on $\Gamma_0(N)$, called ``newforms'', have a simple
characterization in terms of their Fourier coefficients and
we only search for newforms in our calculations.

In this paper we consider Maass forms on $\Gamma_0(p)$ for $p$ prime.
In this case there are four symmetry types:
$(even,+)$, $(even,-)$, $(odd,+)$, and $(odd,-)$.  For $p\le 107$ 
we exhibit the first eigenvalue in each symmetry type, and 
we found the first 1000 newform eigenvalues on~$\Gamma_0(11)$.

\subsection{Maass form $L$-functions}
We associate $f(z)$ to an $L$-function by the Mellin transform.  
If $f$ is even then 
\begin{equation}\label{eqn:mellineven}
(2\pi)^{-s}L(s)G(s)=\int_0^\infty f(iy)y^{s-\frac32}\,\frac{dy}{y},
\end{equation}
where
$$
L(s)=\sum_{n=1}^\infty \frac{a_n}{n^s}
$$
is the associated $L$-function,
and
\begin{equation}\label{eqn:G}
G(s)=\int_0^\infty K_{iR}(y)y^{s}\frac{dy}{y}=2^{s-2}
\Gamma \left(\frac{s+iR}2 \right) \Gamma \left(\frac{s-iR}2 \right).
\end{equation}
Note that the series and integrals converge if $\Re(s)$ is
sufficiently large.

If $f$ is odd then 
\begin{equation}\label{eqn:mellinodd}
(2\pi)^{-s}L(s)G(s+1)=\int_0^\infty 
\frac{\partial f}{\partial x}(iy)y^{s-\frac12 }\,\frac{dy}{y},
\end{equation}
with $L(s)$ and $G(s)$ as before.

The key properties of $L(s)$ are summarized in the following

\begin{proposition}  $L(s)$ continues to an entire function which satisfies the 
functional equation
\begin{eqnarray}\label{eqn:functionalequation}
\xi(s)&=&N^{-{\frac14}+\frac{s}{2}}(2\pi)^{-s}L(s)G(s+a)\cr
&=&\pm (-1)^{a+1}\xi(1-s),
\end{eqnarray}
where $a=0$ if $f$ is even and $a=1$ if $f$ is odd,
the $\pm$ is determined by the eigenvalue of $f$ under the Fricke
involution and whether $f$ is even or odd,
and $G(s)$ is given in~(\ref{eqn:G}).
\end{proposition}

\begin{proof}
Let $h$ equal $f$ or $\partial f/\partial x$, depending on whether 
$f$ is even or odd, respectively.  Combining 
(\ref{eqn:mellineven}) or (\ref{eqn:mellinodd}) with (\ref{eqn:fricke})
we have
\begin{equation}\label{eqn:hecketrick}
(2\pi)^{-s}L(s)G(s+a)=\int_\Delta^\infty h(iy)y^{s-\frac{3}{2}+a}dy
\mp
(-1)^{a} N^{\frac12 -s} \int_{\frac{1}{N \Delta}}^\infty
h(iy)y^{-s-\frac12 +a}dy
\end{equation}
for any $\Delta >0 $.   By (\ref{eqn:kdecay}) the above integrals converge for any
value of~$s$, which proves the analytic continuation.  
Set $\Delta=1/\sqrt{N}$
to see that the right side satisfies the 
functional equation.
\end{proof} 

\section{Locating Maass forms}\label{sec:thealgorithm}

The methods we use to locate Maass forms are similar to previous
methods which have been used, such as that of Hejhal~\cite{Hej}.
In this paper we only consider the case of $\Gamma_0(p)$, 
$p$~prime, and furthermore we use the Fricke involution to 
separate symmetry types, so this means that we do not have to
model vanishing at the cusps.  Thus, except for the larger
number of generators we can in principle use the exact same
methods that have been used for Hecke triangle groups.

The following is a summary of the methods used to locate an individual
Maass form.  
Given generators $\{g_j\}$ of $\Gamma$, we produce an overdetermined
system of linear equations which uses a truncation of the
Fourier expansion of $f$
$$
\tilde{f}(z)=\sqrt{y}\sum_{|n|\leq M,n\not=0} a_n  K_{iR}(2 \pi |n| y)\exp(2 \pi i n x),
$$ 
where the $a_i$ are complex unknowns. Note that we assume $a_0=0$
(which excludes all Eisenstein series when we have only one cusp) and
also normalize one of the coefficients (usually $a_1$) to equal~$1$.
Also note that this of course depends on $R$ (or equivalently on
$\lambda$).  For most of the examples in this paper, we choose~$M$
so that the error caused by the truncation is around $10^{-8}$
for points in the fundamental domain of~$\Gamma$.

We treat the $\{a_n\}$ as $4M-2$
real unknowns. Next we choose $N$
points $z_i$ (where $N>4M$)  on a horizontal line in $\H$. These
points are mapped by the generators to points $g_j z_i=z_i^*$ higher up in
$\H$. If $f$ is a Maass form on $\H/\Gamma$ then $f(z_i)=f(z_i^*)$ (or
more generally $f(z_i)=\chi(g_j)f(z_i^*)$ where $\chi$ is a 
character). The $N$ equations 
$$ 
\tilde{f}(z_i)=\tilde{f}(z_i^*)
$$
constitute an overdetermined system $Ax=b$ in $4M-2$ unknowns.
If $R$ is an eigenvalue of a Maass form on $\Gamma$, then this
system should be consistent to within the error caused by the truncation.
Note: if one wishes to search only for newforms, then it is necessary to
include additional equations to specify this.  See~\cite{R}.

Next we determine
the least square solution $\tilde x$ (using $QR$-factorization) to
this system of equations. We then use the norm of the error,
$||A\tilde{x}-b||_2$, as a measure of how close $\lambda$ is to an
eigenvalue. If $\lambda$ is really an eigenvalue, then
$||A\tilde{x}-b||_2$ should be roughly the size of the truncation
error. In our initial calculations in cases where earlier data were
available we found this to be true. We also found that away from
eigenvalues (i.e. if we choose $R$ randomly) the error
is generally of size~$1$, independent of the size of the truncation error.
We take $||A\tilde{x}-b||_2$ to be a measure of distance
between $R$ and a ``true'' eigenvalue for~$\Gamma$, and we have
found this measure to vary smoothly and to be consistent
with various other checks, which we describe below.  Thus, it seems reasonable to
say that these functions are a factor of $10^8$ closer to
being invariant under~$\Gamma$ then a randomly chosen function.

There are a number of error checks. If $\Gamma$ is arithmetic
then the Fourier coefficients will be multiplicative, and we
find that the above method produces functions whose coefficients
are multiplicative to better than  the truncation error.
These can be viewed as independent 1~in~$10^8$ error checks,
which render it very likely that the functions produced are
indeed Maass forms.  In other words, the possibility of a ``false alarm''
is extremely small, and we have high confidence that the program
is finding the Maass forms for~$\Gamma$.

For nonarithmetic groups there are no Hecke relations, but there
are other persuasive checks.   We start with a general Fourier
expansion with complex coefficients.  For the Maass forms we find,
the functions are real to very high accuracy (to an even higher
accuracy then the truncation error).  In general, when we are far from
an eigenvalue, the system of equations is far from consistent
and the approximate solutions are far from real.

A final check is the size of the Fourier coefficients.  For
arithmetic~$\Gamma$ all the coefficients we have found fit the
Ramanujan-Peterson conjecture~$|a_p|<2$.  For nonarithmetic
groups that bound is not true in general, but it is
still conjectured that~$a_n\ll n^\epsilon$.  And we do in
fact find that if $R$ is close to an eigenvalue then
the $a_n$ from the least-squares solution are much smaller
(and not growing as a function of $n$)
than those from random~$R$.

Our programs were implemented in Mathematica.

\section{Low eigenvalues for $\Gamma_0(p)$}\label{sec:loweigenvalues}

In Table~\ref{tab:lowest} we present the first eigenvalue in
each of the four symmetry types for $\Gamma_0(p)$, $p$~prime.

\begin{table}[h!tb]
\small
\begin{tabular}{r|l|l|l|l}
\multicolumn{1}{c|}{$p$} & \multicolumn{1}{c|}{$(even,+)$} & 
\multicolumn{1}{c|}{$(even,-)$} & \multicolumn{1}{c|}{$(odd,-)$} & 
\multicolumn{1}{c}{$(odd,+)$} \\
\hline\hline
2 & 8.9228764869917 
& 12.092994875078
& 7.220871975958 
&  5.4173348068447\\

3 & 5.0987419087295
& 8.7782823935545
& 6.1205755330872
& 4.3880535632221\\

5 & 4.1324042150632
& 5.436180461416
& 4.897235015733
& 3.028376293066\\

7 & 3.454226503571
& 4.8280076684720
& 4.119009292925
& 1.924644305111\\

11 & 2.4835910595550
& 4.018069188221
& 2.96820576382
& 2.033090993855\\
\hline
13 & 2.025284395696
& 3.701627575242
& 2.8308066514473
& 0.97081541696\\

17 & 1.849687906031
& 3.169382380088
& 1.967986359638
& 1.441428545022\\

19 & 1.32979889046
& 3.0371960205
& 2.297006359074
& 1.09199155992\\

23 & 1.57958924015
& 2.61095996203
& 1.393337141483
& 1.5061266371\\

29 & 1.01726655080
& 2.35848525400
& 1.4875542669
& 1.21206072002\\
\hline
31 & 0.78935617774
& 2.3681029381
& 1.68678370214
& 1.06284037124\\

37 & 1.22324304054
& 2.07459336618
& 1.7095550812
& 0.6423059582\\

41 & 0.66572483212
& 2.00647646730
& 1.204395742494
& 0.86739746584\\

43 & 1.10814196343
& 1.80282682958
& 1.33429136841
& 0.65545238186\\

47 & 1.11157408467
& 1.6028012074
& 0.5854521430
& 1.2635566239\\
\hline
53 & 1.0158404282
& 1.49379571497
& 1.10150481917
& 0.8039894596 \\

59 & 0.59582968816
& 1.7096938151
& 0.5985969867
& 1.0391204647\\

61 & 0.67704450723
& 1.7475722491
& 1.40716257846
& 0.41806231115\\

67 & 0.84247476
& 1.3611564746
& 1.0195601623
& 0.6775902092\\

71 & 0.3574504906
& 1.579666665
& 0.5829071706
& 1.02176122015\\
\hline
73 & 0.654548913
& 1.402091810
& 1.3304765833
& 0.5178875052\\

79 & 0.5517748428
& 1.5496128514
& 0.9963893608
& 0.7490008694\\

83 & 0.817132725
& 1.0975181181
& 0.6402875784
& 0.95188877503\\

89 & 0.489435980
& 1.4372631258
& 0.860333981
& 0.688804233\\

91 & 0.6038860327
& 1.206246429
& 1.132502360
& 0.416931671\\
\hline
101 & 0.453759608
& 1.30563134
& 0.546229605
& 0.77990698\\

103 & 0.685190643
& 1.08157342
& 0.724051309
& 0.56540244\\

107 & 0.840011226
& 0.90440769
& 0.581677094
& 0.90574018\\
\end{tabular}

\caption{\sf \label{tab:lowest}
First eigenvalue, $R$, in each of the four symmetry types of
$\Gamma_0(p)$, for prime $p\le 107$.  Values shown are truncations
of the actual values, and all digits are believed to be correct.}
\end{table}

Table~\ref{tab:lowest} shows a few trends.  For example,
the first $(even,+)$ eigenvalue tends to be smaller
than the first $(even,-)$ eigenvalue.  This is not
surprising because the number of nodal domains is an
increasing function of the eigenvalue, and the symmetry
relations force various constraints on the nodal curves.
This also has an effect on the lower order terms in
Weyl's law, which can be seen in Table~\ref{tab:first30}.

\begin{figure}[htp]
\begin{center}
\scalebox{1.5}[1.5]{\includegraphics{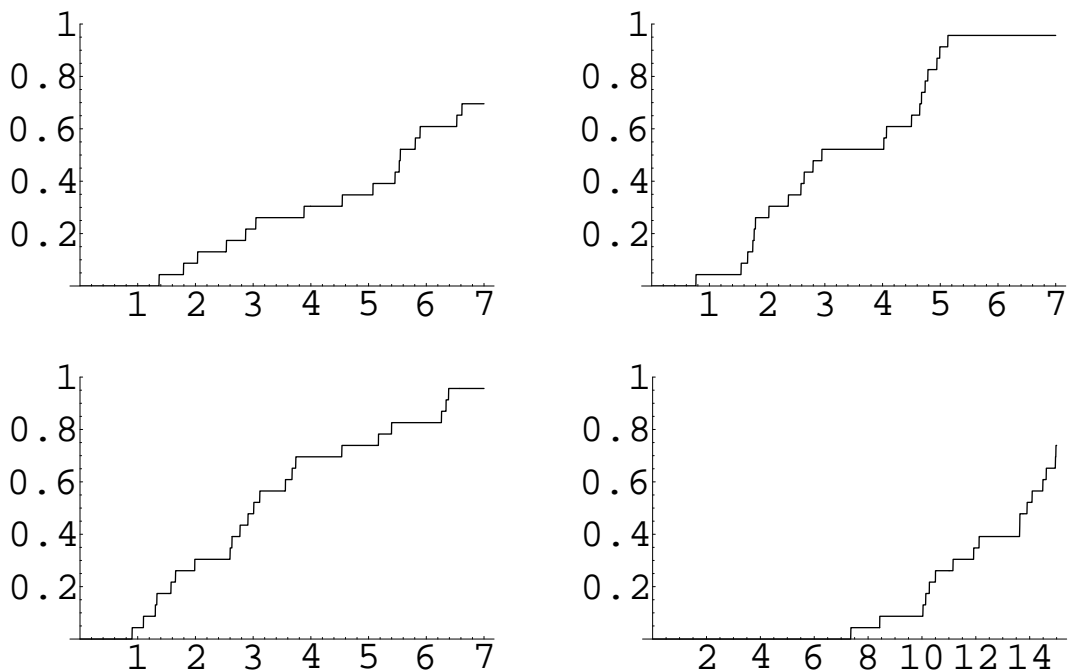}}
\caption{\sf\label{fig:rescaledfirsteig}
Cumulative distributions of the rescaled first eigenvalue
in each of the four symmetry types for $\Gamma_0(p)$.
Clockwise from top right we have even plus, even minus, odd minus, and
odd plus.
} 
\end{center}
\end{figure}

To further study the low eigenvalues, we considered the 
normalized eigenvalues $\tilde\lambda = \frac{p+1}{12}R^2$, 
where $\frac{p+1}{12}$ is the (reciprocal of the)
coefficient of the leading
term in Weyl's law, so that $\tilde\lambda_n\sim n$ for the
$n$th eigenvalue.  
Figure~\ref{fig:rescaledfirsteig} shows the cumulative distribution
functions of the first $\tilde\lambda$ in each symmetry type.
If each symmetry type was equally likely to have the lowest
eigenvalue then each plot would correspond to a p.d.f with
mean~$4$.  This is not the case, for the reasons described 
in the previous paragraph.  A more sophisticated normalization,
involving the lower order terms in Weyl's law, would be 
needed in order to reveal any underlying structure.

\section{The first 1000 Maass forms on $\Gamma_0(11)$}\label{sec:gamma011}

We applied the method described in the previous section to
find all Maass newforms on $\Gamma_0(11)$ with $R<36.5$,
resulting in a total of 1054 eigenvalues.
The first 30 eigenvalues in each symmetry type are given
in Table~\ref{tab:first30}.

\begin{table}[h!tb]
\small
\begin{tabular}{l|l|l|l}
\multicolumn{1}{c|}{$(even,+)$} &
\multicolumn{1}{c|}{$(even,-)$} & \multicolumn{1}{c|}{$(odd,-)$} &
\multicolumn{1}{c}{$(odd,+)$} \\
\hline\hline
2.48359105931 & 4.01806918817 & 2.03309099399 & 
    2.96820576406\\ 
3.28347243577 & 4.55526078519 & 3.03251283477 & 
    3.67969872541\\ 
4.71167698688 & 5.763876683 & 3.48188847623 & 
    4.52889442543\\ 
4.93319801514 & 6.24170090061 & 4.47950391861 & 
    5.38428160742\\ 
5.26506649404 & 6.50705758603 & 4.73801386035 & 
    5.86853222013\\ 
\hline
5.98954160832 & 7.16536040411 & 4.96681056818 & 
    6.03361084498\\ 
6.53884437369 & 7.73368897382 & 5.66326548842 & 
    6.77549239897\\ 
6.61367515372 & 7.88007383371 & 6.06921835945 & 
    7.02463370057\\ 
7.26644184877 & 8.62006794999 & 6.58162198432 & 
    7.59926619568\\ 
7.57048468847 & 8.99092809554 & 6.83944314323 & 
    7.68091057664\\ 
\hline
8.1420217952 & 9.10172820901 & 7.29898901659 & 
    8.43825975028\\ 
8.2745067789 & 9.22640277057 & 7.48914848604 & 
    8.55580430772\\ 
8.51378467778 & 9.58661569223 & 7.51799881213 & 
    8.57197306999\\ 
8.67868612906 & 10.1232769351 & 8.04940484112 & 
    9.04133195689\\ 
9.34515203418 & 10.3790777119 & 8.16066289206 & 
    9.46659932758\\ 
\hline
9.58209996953 & 10.6466817416 & 8.65213300992 & 
    9.58865618615\\ 
9.67447022214 & 10.9372887384 & 8.92993868939 & 
    10.1617823366\\ 
9.86213636813 & 11.1472353004 & 9.19975228208 & 
    10.2760011943\\ 
9.98570867995 & 11.4784159679 & 9.67408869747 & 
    10.6445164137\\ 
10.6430998724 & 11.785015526 & 9.70588526636 & 
    10.6694194813\\
\hline
10.6895930795 & 11.8109636758 & 9.85493309024 & 10.989904792 \\
11.1502516477 & 12.0601741852 & 10.1854915198 & 11.2612425027 \\
11.1780820772 & 12.6105517277 & 10.3386794422 & 11.2733375812 \\
11.4741161917 & 12.6118028204 & 10.5313292208 & 11.8005163423 \\
11.6836999337 & 12.9273387317 & 10.8165867701 & 11.8348731714 \\
\hline
11.9293065449 & 12.9741004148 & 10.9051400611 & 11.9611740135 \\
12.038113214 & 13.1014471689 & 11.3280653248 & 12.2747697039 \\
12.3049680989 & 13.3464693158 & 11.5395912507 & 12.6217030476 \\
12.5345178216 & 13.9313207377 & 11.6712226728 & 12.8306123388 \\
12.6348846255 & 13.9973467855 & 11.7801194448 & 13.0161442996\\
\end{tabular}

\caption{\label{tab:first30}
First 30 newform eigenvalues in each symmetry type of $\Gamma_0(11)$.
Values shown are truncations of the actual values, and all digits
are believed to be correct.
Note that there are only three oldforms in the range covered be the table:
an even oldform at $R\approx 13.779751$, and
odd oldforms at $R\approx 9.533695$ and $R\approx 12.173008$.
}
\end{table}

\subsection{Weyl's law justification that the list is complete}

It is natural to question whether our list of 1000 eigenvalues is
complete.  It is possible to do a calculation using the trace
formula to justify this assertion, as in~\cite{BSV}.
However, we will use a simpler but less rigorous method
involving Weyl's law.  This approach is well known to
physicists 
and we do not claim that any of the ideas in this section are new.

The method we describe is quite robust and does not require
any {\it a priori} knowledge about the eigenvalue counting
function.  In addition to checking lists of eigenvalues, the
method also works well for zeros of $L$-functions, 
zeros of derivatives of $L$-functions, and no doubt many
other cases.  In some of these cases we know very precise information
about the counting function of the objects being studied, but this
is not necessary for the method.

Let $N(x)=\#\{\lambda_j\ :\ 0<\lambda_j \le x\}$.
The method is based on the assumption that there exists
a nice function $f$, which is of a fairly simple form,
such that
\begin{equation}\label{eqn:weylfit}
r(x):=N(x)-f(x)
\end{equation}
grows very slowly and averages to~$0$ over small intervals.
Weyl's law provides the leading order behavior of~$f(x)$,
and in many cases some lower order terms are known,
but we do not claim that one can \emph{prove} that the
list of eigenvalues is complete only by using
established facts about Weyl's law.

We justify that our list of eigenvalues
is complete in two steps: 
\begin{itemize}
\item demonstrate that the function $f$ exists for our
claimed complete list of eigenvalues.
\item demonstrate that the function $f$ does not exist 
for lists with missing eigenvalues.
\end{itemize}

\begin{figure}[htp] 
\begin{center}
\scalebox{0.75}[0.75]{\includegraphics{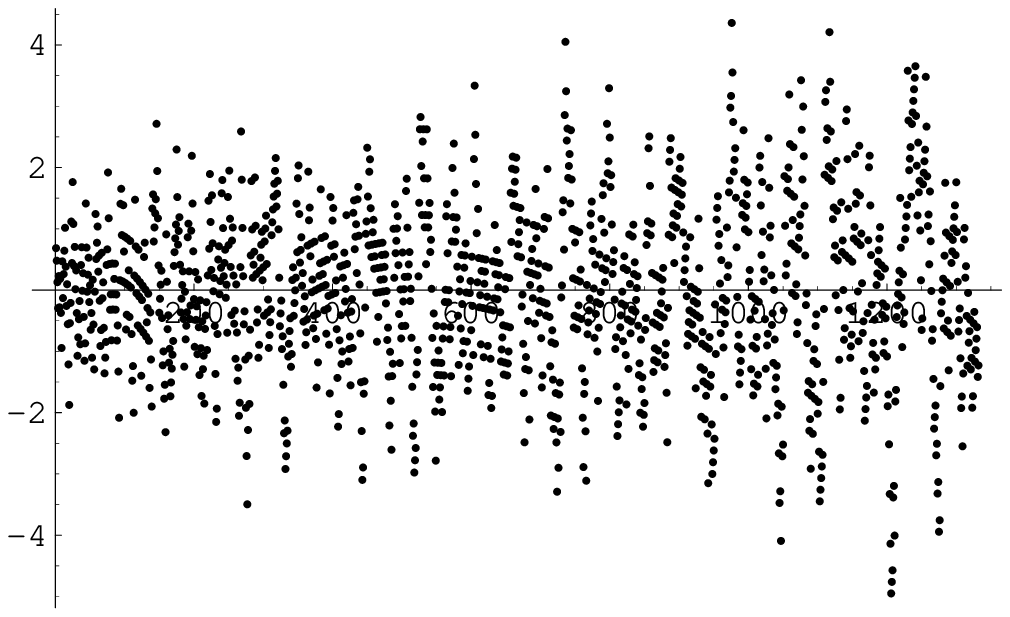}}
\hskip 0.3in
\scalebox{0.75}[0.75]{\includegraphics{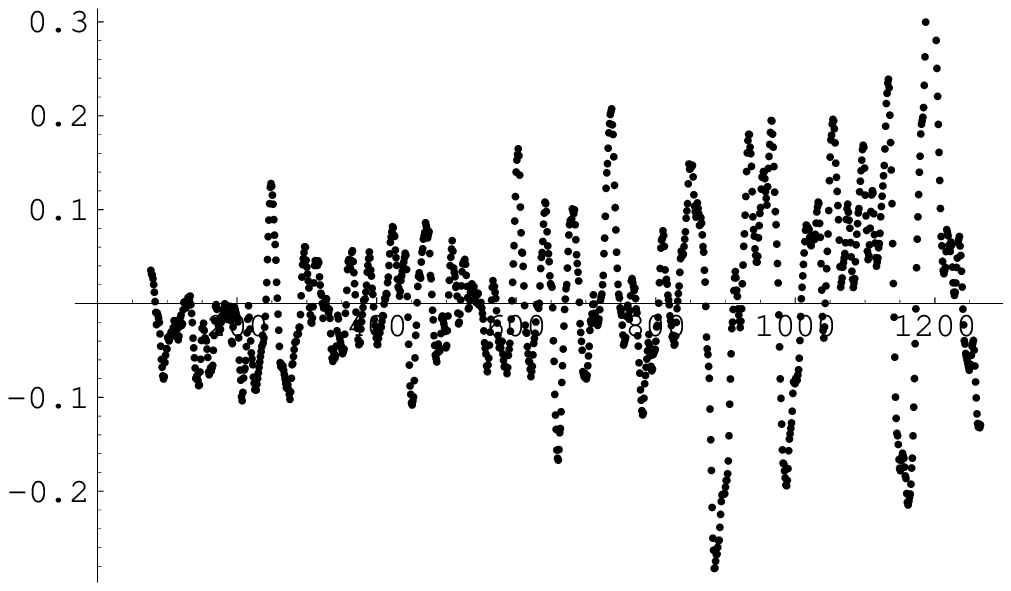}}
\caption{\sf
The remainder function $r(x)$ defined in~\eqref{eqn:weylfit}
where $f(x)$ is chosen as described in the text.
The data is the 1054 newforms on $\Gamma_0(11)$
with $\lambda<1260$.
The plot on the right is the running average over a window
of length 150.
} \label{fig:weylall}
\end{center}
\end{figure}

First we consider the case of all newforms on $\Gamma_0(11)$,
and we assume that the function $f$ has the form 
$f(x)=A x + B \sqrt{x}$.
If we choose $A$ and $B$ to give the best least-squares fit
we find
$A=0.83305$ and $B=-1.51762$.  
This shows good agreement with the known main term
of $A=10/12\approx 0.83333$.
Figure~\ref{fig:weylall} plots the function $r(x)$, as well
as $r(x)$ averaged over a moving window of width~150.

As can be seen, the function $r(x)$ is small, and it is also
small on average. 

\begin{figure}[htp]
\begin{center}
\scalebox{0.75}[0.75]{\includegraphics{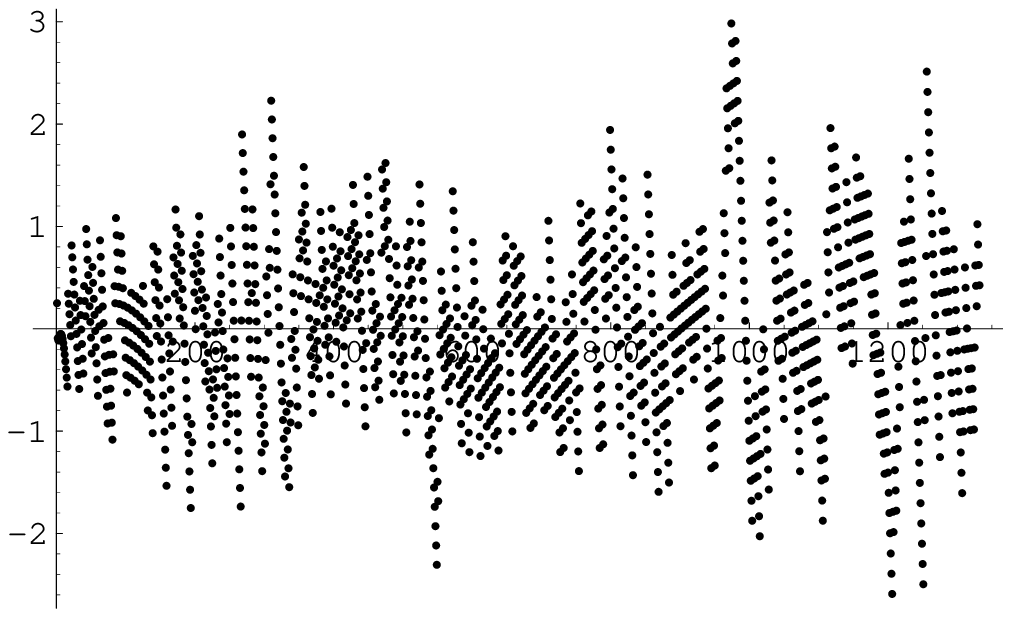}}
\hskip 0.3in
\scalebox{0.75}[0.75]{\includegraphics{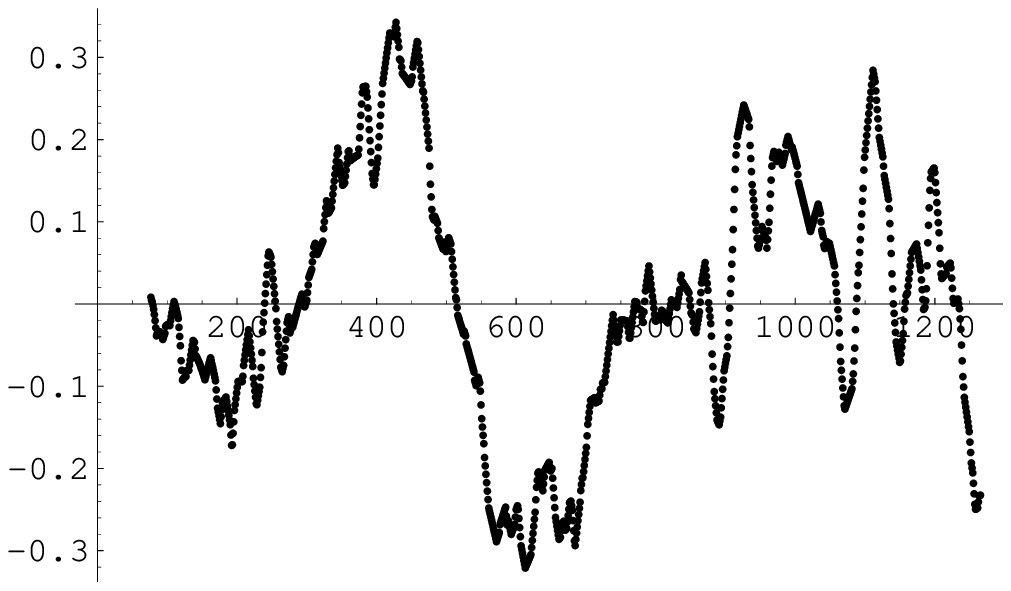}}
\caption{\sf
The remainder function $r(x)$ for the first 250 (even,$-$) 
newform eigenvalues on $\Gamma_0(11)$, where
{\em {one eigenvalue with $R^2\approx 510$ has been intentionally omitted}}.  
The plot on the right is the running average over a window
of length 150.
} \label{fig:weylem}
\end{center}
\end{figure}

We now repeat the calculation on a list of eigenvalues which 
is known to be incomplete, to demonstrate how the 
method detects the missing value.  In Figure~\ref{fig:weylem}
we consider the eigenvalues in the (even, $-$) space,
for which we have removed one eigenvalue~$\lambda\approx 510$.
Since we are in the (even, $-$) space, there are only 1/4
as many eigenvalues, so it should actually be easier to 
fit a main term function~$f(x)$.  We want to make the point that
is impossible to find functions $f$ and $r$ for this list
of data.  We choose $f$ of the form
\begin{equation}
f(x)=A x + B \sqrt{x}\log x + C \sqrt{x} + D x^{\frac14}
+ E \log x + F,
\end{equation}
and then choose $A,B,C,D,E,F$ to give the best least-squares fit.
The resulting $r(x)$ and $r(x)$ averaged over a window of
size~150 is shown in Figure~\ref{fig:weylem}.

As can be seen, $r(x)$ does not average to zero over small intervals,
and in fact we can almost read off the missing value 
$\lambda\approx 510$ from the averaged graph.
The fact that we fit a function with many free parameters is
meant to suggest the robustness of the method.
Unfortunately, when fitting such a general function to such a small 
amount of data, we can no longer use the coefficient $A$ of the
main term to check the known main term in Weyl's law.


\subsection{Shortcomings of the method}

There are some obvious shortcomings of the method of fitting
Weyl's law to check the data.  

First, the method cannot detect
missing eigenvalues at the beginning of the list.  Such
missing term are absorbed by the constant term in our
function.  

Second, if a large number of eigenvalues are
missing, and the missing eigenvalues are very regularly spaced,
then this may not be detected.  Indeed, that is exactly
the situation at hand.  We are claiming to demonstrate that
our list of {\em{newform}} eigenvalues is complete.  But
if we were mistakenly claiming that we had the list
of {\em{all}} eigenvalues, that error would not be immediately caught 
by our fit to ``Weyl's law,'' for the missing eigenvalues, 
i.e. the oldforms, are regularly spaced.  Indeed, they
have their own Weyl's law.   In the case here, we know the 
coefficient of the
main term in the counting function 
for newforms, and our fit function shows good agreement with that.
However, we cannot rule out the possibility that we are
systematically missing every 100th eigenvalue, for example.

\subsection{Statistics of the eigenvalues}

It is generally believed 
that the newform eigenvalues in each symmetry class are as uncorrelated
as possible (have Poisson statistics),  the 
newform eigenvalues in the different symmetry classes
are uncorrelated,  and there is no correlation between newform
eigenvalues and oldform eigenvalues.  For the case of 
$SL(2,\Z)$ these conjectures are strongly supported by numerical
evidence~\cite{Hej}.  We now check if our
data supports these conjectures.
Only the case of the spacing of oldforms within the list of newforms
can be considered new here.

\begin{figure}[htp]
\begin{center}
\scalebox{0.7}[0.7]{\includegraphics{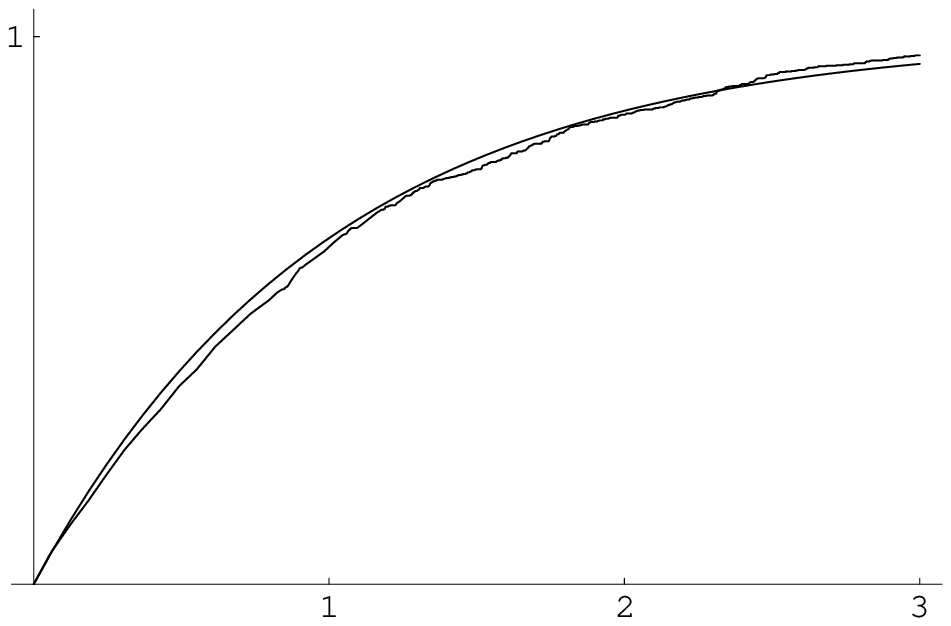}}
\hskip 0.3in
\scalebox{0.7}[0.7]{\includegraphics{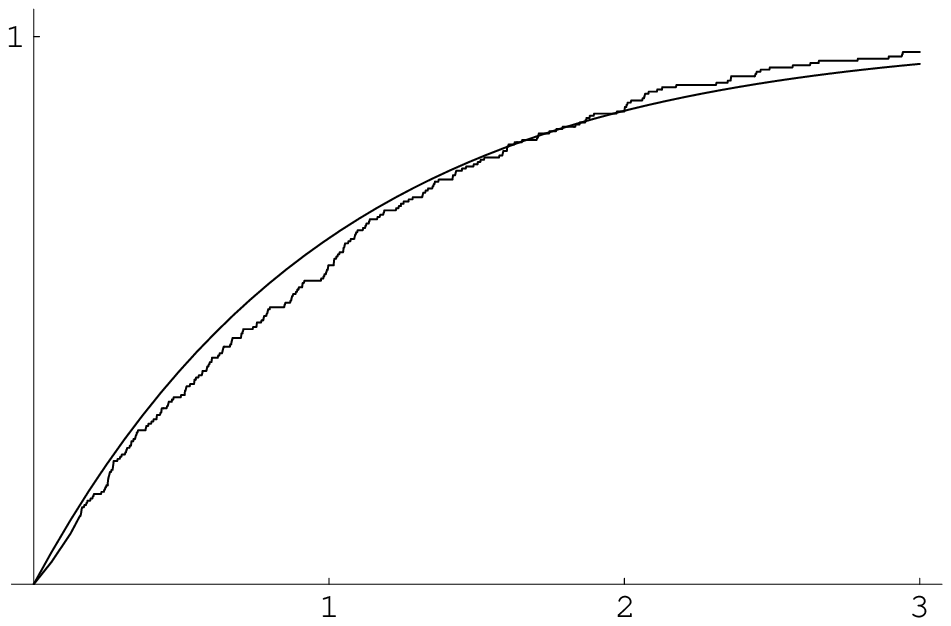}}
\caption{\sf
Cumulative distribution of normalized nearest neighbor spacing 
of newforms on $\Gamma_0(11)$.  The smooth curve is the 
cumulative distribution for Poisson spacing.
The plot on the left is for all newforms and the plot on the right is
for (even,$-$) forms.
} \label{fig:eignn}
\end{center}
\end{figure}

In Figure~\ref{fig:eignn} we consider the nearest neighbor spacing of
the eigenvalues on $\Gamma_0(11)$, both for the full spectrum and
for one particular symmetry type.  Good agreement with Poisson 
statistics is found.

\begin{figure}[htp]
\begin{center}
\scalebox{0.70}[0.70]{\includegraphics{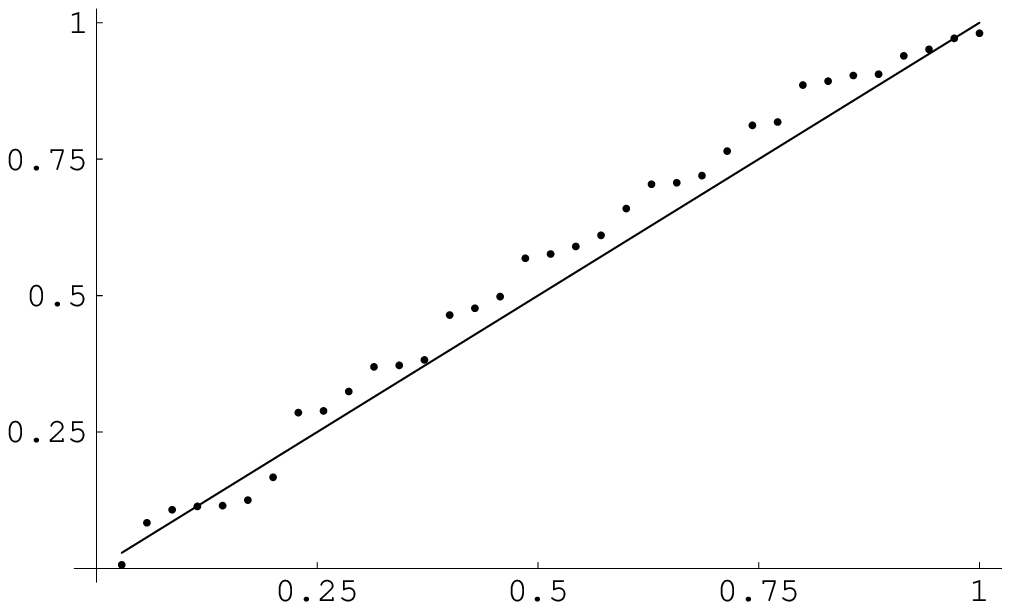}}
\hskip 0.3in
\scalebox{0.70}[0.70]{\includegraphics{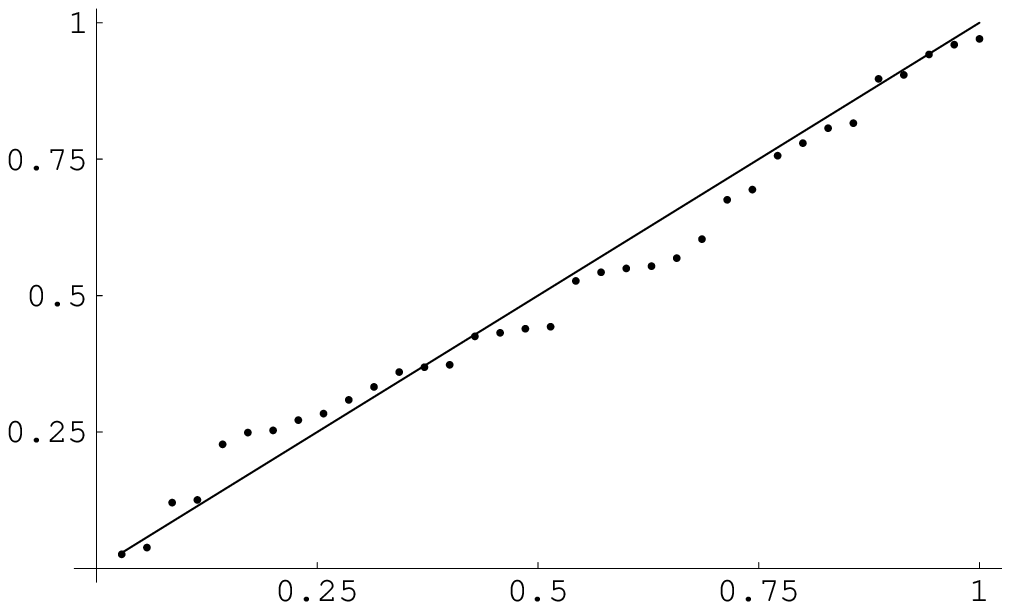}}
\caption{\sf
Ordered normalized location of odd oldforms within neighboring newforms, 
for the 35 oldforms with $R<36.5$.   
Straight line is comparison to random (Poisson) spacing.  
The plot on the left is for odd$-$,
on the right for odd$+$.  
} \label{fig:oldodd}
\end{center}
\end{figure}

Next we consider the location of the oldforms in among the newforms.
If there is no correlation then the oldforms should be distributed
uniformly between the two neighboring newforms.  
Figure~\ref{fig:oldodd} shows the (scaled) location of the
oldforms between the neighboring newforms, giving a good
agreement with the expected Poisson distribution.

\section{Example Maass $L$-functions}\label{sec:examples}

We give example plots of Maass $L$-functions and discuss some
general features, and in the next section we consider the statistics
of zeros and critical values.

\begin{figure}[htp]
\begin{center}
\scalebox{0.74}[0.74]{\includegraphics{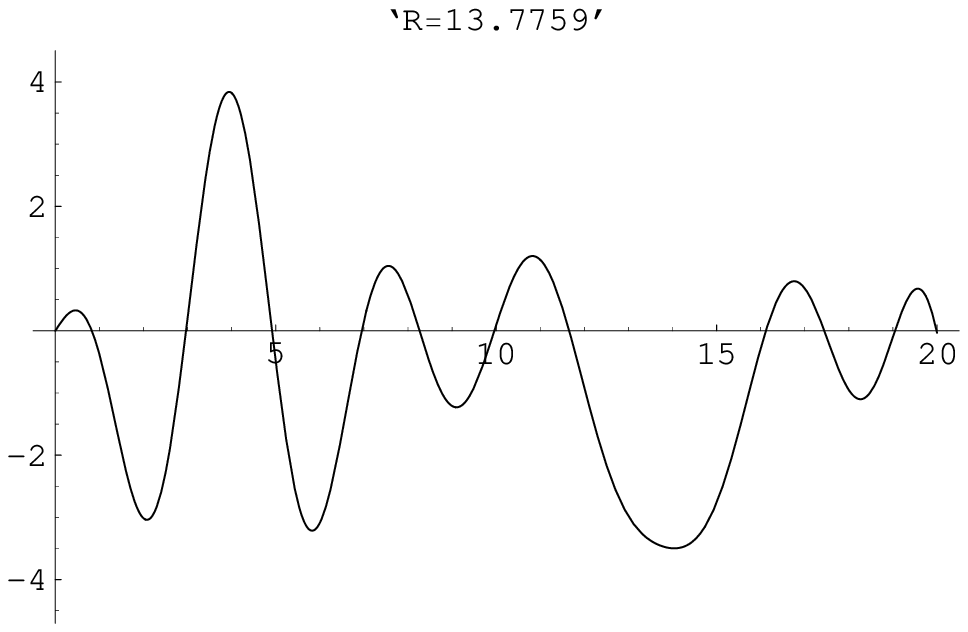}}
\hskip 0.3in
\scalebox{0.74}[0.74]{\includegraphics{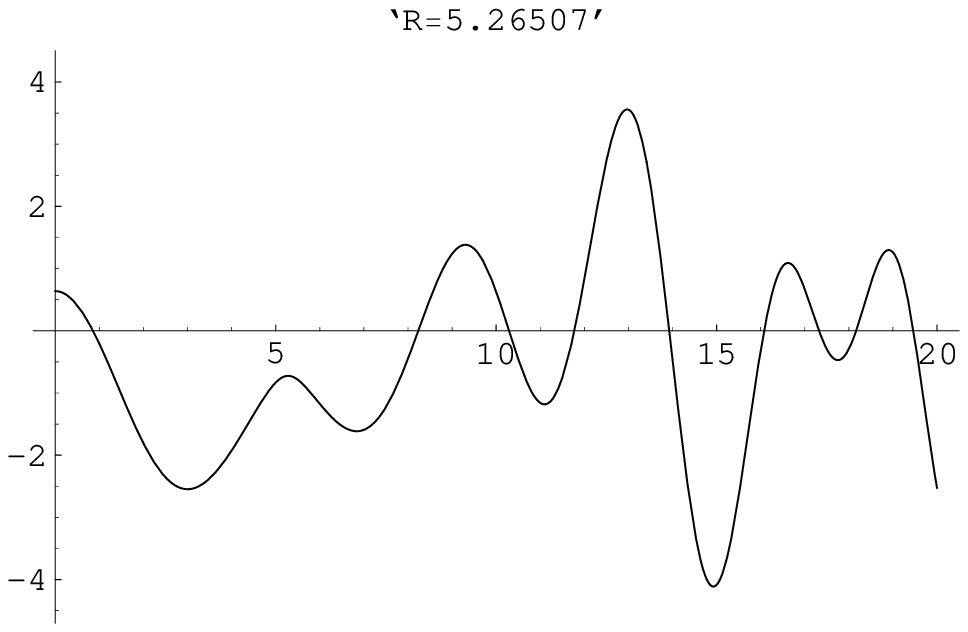}}
\vskip 0.3in
\scalebox{0.74}[0.74]{\includegraphics{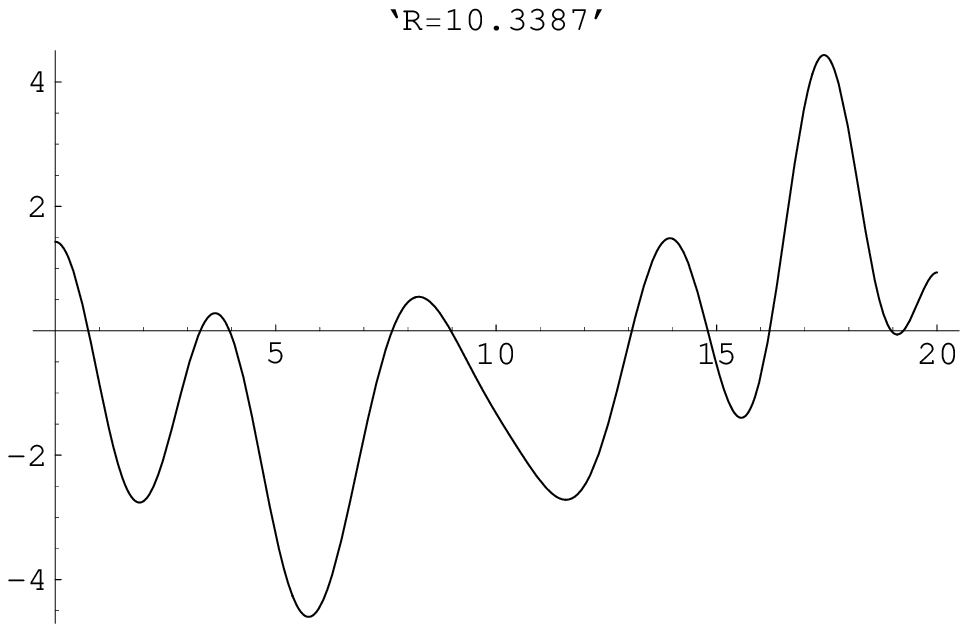}}
\hskip 0.3in
\scalebox{0.74}[0.74]{\includegraphics{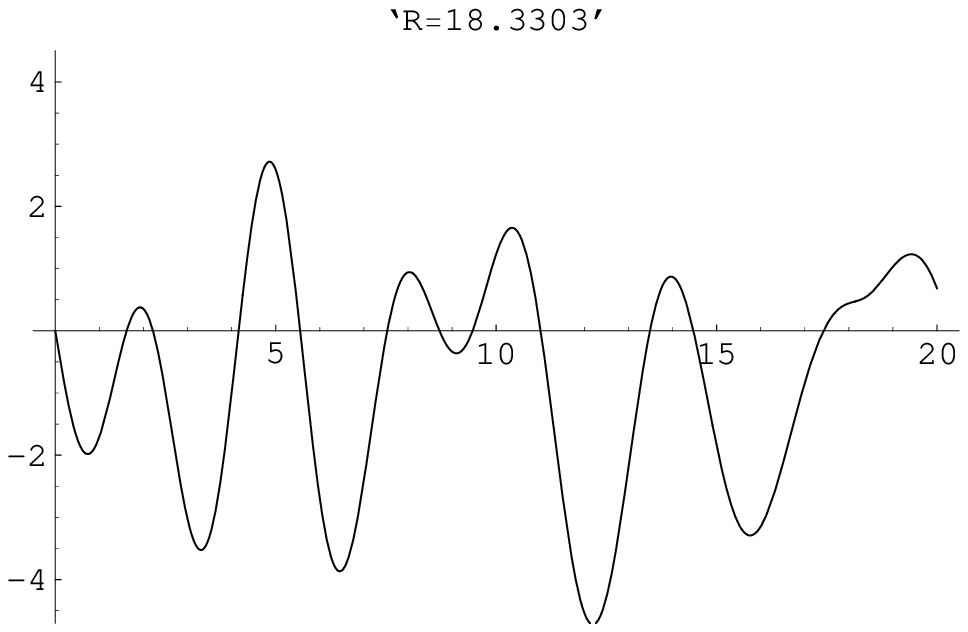}}
\caption{\sf 
Plots of (the real version of) $L(\frac12+it)$, one for each of the four symmetry types
of Maass forms on $\Gamma_0(11)$.
Clockwise from top right we have odd plus, even plus, even minus, and
odd minus.  
Each plot is labeled with the``eigenvalue'' $R$ of the Maass form.
} \label{fig:fourplots}
\end{center}
\end{figure}

By the functional equation, $\xi(\frac12+it)$ is real for real~$t$,
so there is a simple rescaling (involving a ratio of $\Gamma$-functions)
to produce a function with the same modulus as $L(\frac12+it)$ which
is real for real~$t$.  It is that real function of the real variable~$t$,
which we still call ``the $L$-function'', which is shown in
Figure~\ref{fig:fourplots}.  There are four classes of $L$-functions,
depending on whether the function $f$ is even or odd as a function
of~$x$, and whether $f$ has eigenvalue $+1$ or $-1$ under the
Fricke involution.  One example of each combination is shown
in Figure~\ref{fig:fourplots}.

The first thing to notice is that these plots make it appear that the
Riemann hypothesis is false for these $L$-functions.  For example,
in the upper-right plot, the negative local maximum near $t=5$ 
indicates a zero off the critical line.  However, that zero off
the line is a trivial zero, which has imaginary part $iR$, 
so this is not a counterexample to the Riemann Hypothesis.  

The effect of the trivial zeros is more noticeable in the even
case, where they are closer to the critical line.
This illustrates Str\"ombergsson's~\cite{S}
observation of a gap in critical zeros near $t=R$. 

The average spacing of zeros is a function of both $R$ and $t$,
the precise dependence can be deduced from the functional equation
and the argument principle.  For our purposes it suffices to note that
for $|t|<R$ the zero spacing is primarily a function of $R$, increasing
with~$R$, except for a slight decrease in density near $|t|=R$.

The number $\varepsilon:=\pm (-1)^{a+1}$ from~(\ref{eqn:functionalequation}) 
is called ``the sign of the functional equation.''
If $\varepsilon=+1$ then $L(\frac12)\ge 0$, and we call
$L(\frac12)$ the ``critical value''.
If $\varepsilon=-1$ then $L(\frac12)=0$ and $L'(\frac12)$ is
the critical value.

\section{Checking the random matrix predictions}\label{sec:data}

Random Matrix Theory has become a fundamental tool
for understanding the zeros of $L$-functions. 
Montgomery \cite{Mon} showed that
(in a limited range) the two-point correlations between the 
non-trivial zeros of the
Riemann $\zeta$-function
is the same as that of the
eigenvalues of large random unitary matrices, and he conjectured
that the correlation we in fact equal.
There is extensive
numerical evidence \cite{Odl} in support of this conjecture,
and in particular many statistics of the zeros of the $\zeta$-function, 
such as the distribution of nearest neighbor spacing,
are believed to be the same as that of large random unitary matrices.
Katz and Sarnak \cite{KSa} introduced the idea of
studying zero distributions within families of $L$-functions 
and have conjectured that these coincide with
the eigenvalue distributions of the classical compact groups:
unitary, symplectic, and orthogonal.  For these groups the bulk
of the eigenvalue distributions are the same, it is the 
eigenvalues near~$1$ which have different distributions in
each case.  In terms of $L$-functions this corresponds to zeros
near the critical point, so our concern is with the low lying zeros.

The family of Maass $L$-functions is an Orthogonal family 
(see~\cite{CFKRS}), which we further break into $O^+=SO(even)$ and 
$O^-=SO(odd)$,
depending on whether the $L$-function has sign $+1$ or $-1$ in
its functional equation.  In the latter case $L(\frac12)=0$,
which corresponds to the fact that odd orthogonal matrices have~$1$
as an eigenvalue.

\begin{figure}[htp]
\begin{center}
\scalebox{0.72}[0.72]{\includegraphics{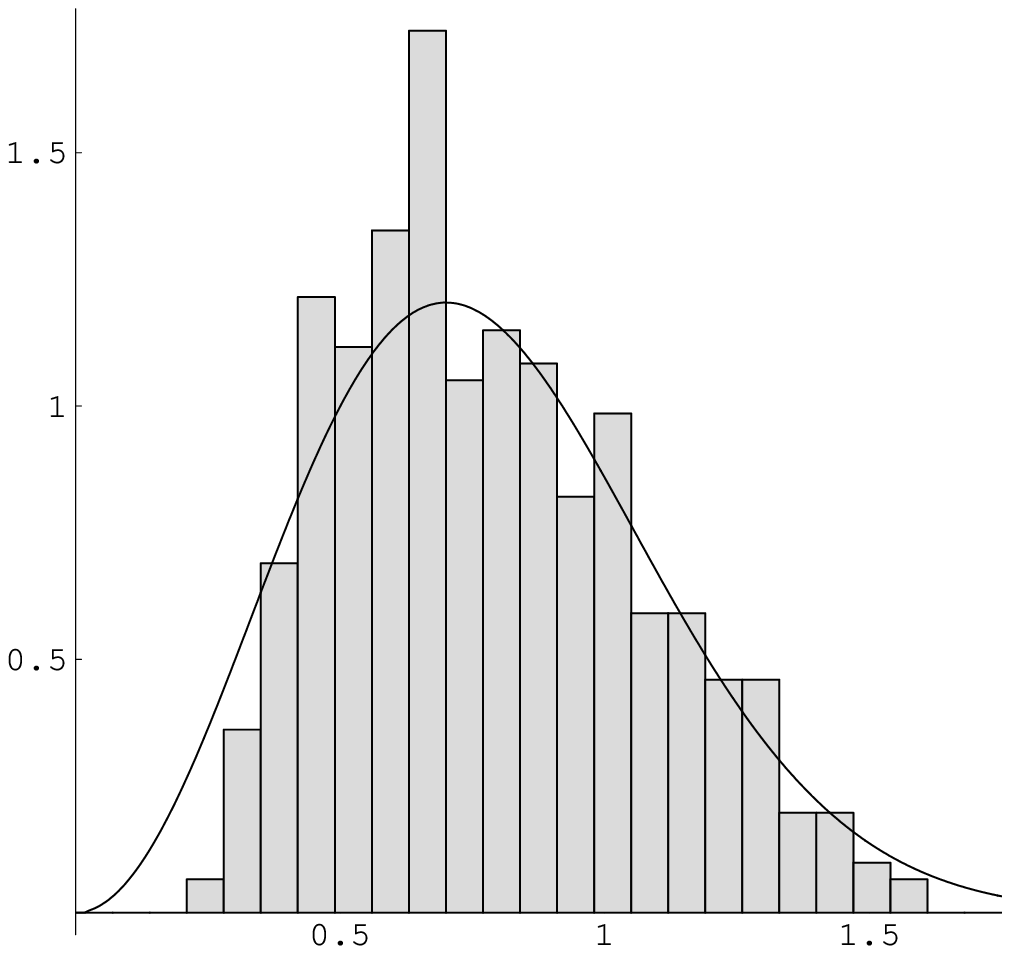}}
\hskip 0.3in
\scalebox{0.72}[0.72]{\includegraphics{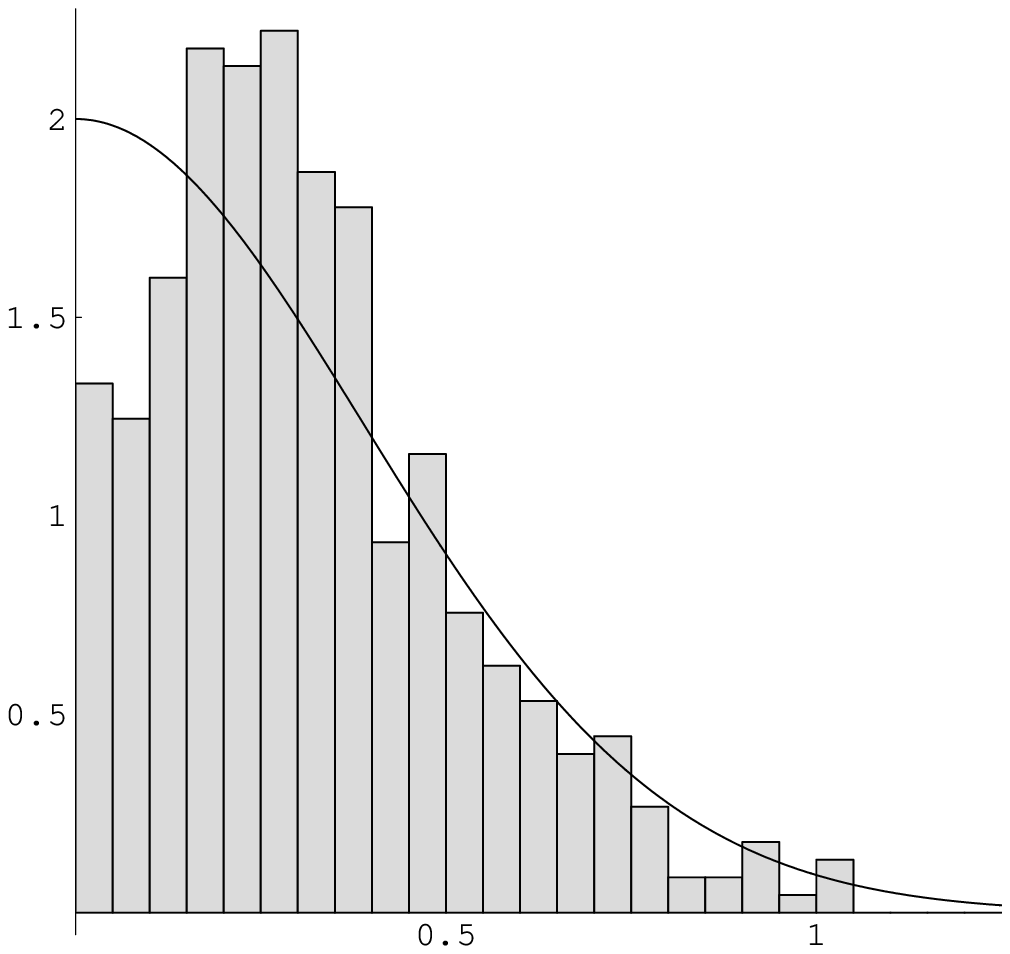}}
\caption{\sf
Density of first zero above the critical point for
435 $L$-functions with $\varepsilon=-1$ on the left, and
450 with $\varepsilon=+1$ on the right. 
The continuous curves are the densities for the corresponding 
eigenvalues for $O^{-}$ and $O^+$, respectively.  The $L$-function
zeros have been rescaled to have the same averages as the
corresponding eigenvalues.
} \label{fig:firstzeros}
\end{center}
\end{figure}

In Figure~\ref{fig:firstzeros} we consider the distribution of the 
height of the first zero above the critical point, comparing to
the analogous quantity for eigenvalues of orthogonal matrices. 
In Figure~\ref{fig:firstzeros} we omit the Maass forms 
with~$R<15$, because
the small values of $R$ cause anomalous behavior 
of the first few zeros.

\begin{figure}[htp]
\begin{center}
\scalebox{1.0}[1.0]{\includegraphics{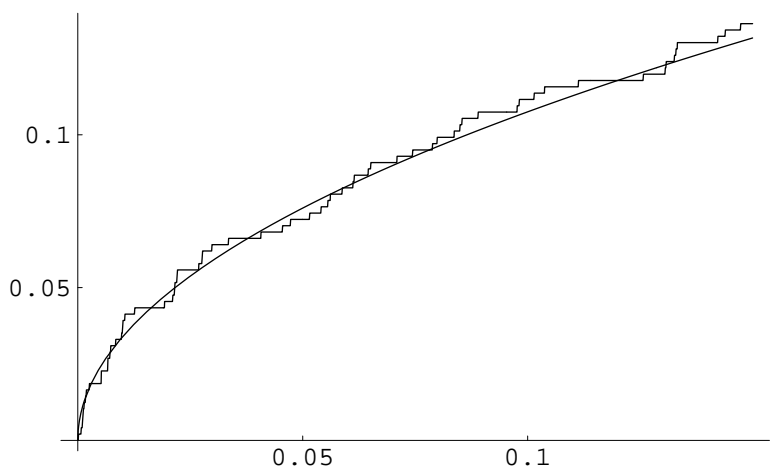}}
\caption{\sf \label{fig:critical1}
Cumulative distribution of critical values for $\varepsilon=+1$ 
$L$-functions. A total of 484 $L$-functions were considered,
of which 63 feature in this plot.  The continuous curve is $0.34\sqrt{x}$.
} 
\end{center}
\end{figure}

Next we show the cumulative distribution of the critical values, fit to
the prediction from random matrix theory.
In the case of even functional equation, corresponding to
$O^+$ matrices,
the density of critical values scales as $c/\sqrt{x}$ for small~$x$.
This is quite difficult to see in a histogram, so instead we consider
the cumulative distribution and compare to the predicted $c'\sqrt{x}$,
where we choose $c'$ by fitting to the data.
Figure~\ref{fig:critical1} shows very good agreement.

In Figure~\ref{fig:critical1} we omit the first 30 Maass forms of
each symmetry type, corresponding to approximately~$R<13$, because
the small values of $R$ cause anomalous behavior near the critical
point.

In the case of odd functional equation, $\varepsilon=-1$, the
density of small values is predicted to scale as $c x^{3/2}$.  In particular,
there should be very few small values.  Indeed we find
very few cases where $L'(\frac12)$ is small, but our statistics
are too low to make a meaningful comparison with the 
random matrix prediction.

\begin{figure}[htp]
\begin{center}
\scalebox{0.8}[0.8]{\includegraphics{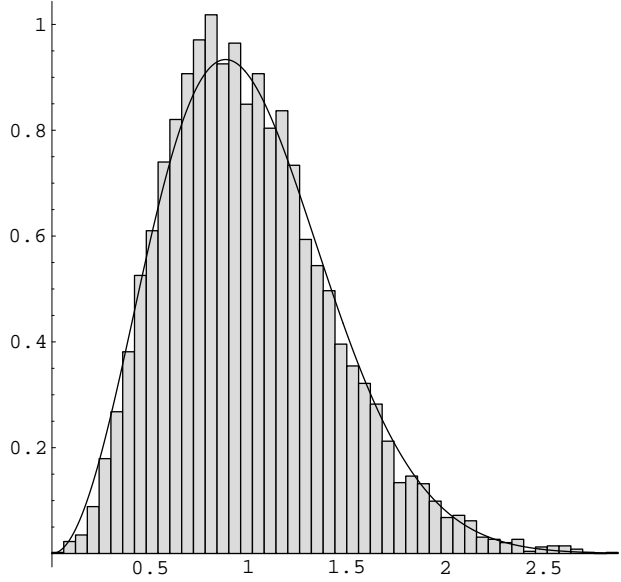}}
\caption{\sf
Nearest neighbor gaps, normalized to remove the dependence
on $R$ and~$t$.
A total of
8086 neighbor gaps from approximately 650 $L$-functions.
The continuous curve is the normalized nearest neighbor spacing
of eigenvalues of large random unitary matrices.
} \label{fig:nearestneighbor}
\end{center}
\end{figure}

Finally, we show the distribution of nearest neighbor spacings.
In Figure~\ref{fig:nearestneighbor} we consider Maass forms with
$R>21.5$, and we considered the neighbor gaps between zeros starting
with the 5th zero and going up to $t=20$.  This
was designed to avoid the anomalous behavior near the critical
point and near~$t=R$.  In addition we rescale the zeros to remove the
trend of decreasing spacing with increasing~$R$ and~$t$. 
The fit to the random matrix prediction in Figure~\ref{fig:nearestneighbor}
is quite good.  Indeed this may be somewhat surprising because the
zeros considered there lie in the range $|t|<R$, which one might
expect to be different from the $t\to\infty$ case.  
Perhaps this illustrates the universality of the random matrix
behavior.

\section{Computing $L(s)$}\label{sec:computing}

Our interest is in the low-lying zeros of $L(s)$, so we only require
an efficient method of evaluating $L(\sigma+it)$ for small~$t$.

Substituting (\ref{eqn:fourier}) in to (\ref{eqn:hecketrick}),
and then switching the order of
summation and integration, we have
\begin{eqnarray}
(2\pi)^{-s}L(s)G(s+a)&=&(2\pi)^{(2a -1)(s)}\sum_{n=1}^\infty
a_nn^{(2a -1)(s)} \int_{2\pi n\Delta}^\infty
K_{iR}(y)y^{s+(a -1)}dy \cr
&&\pm (-1)^{a+1} N^{\frac12-s}(2\pi)^{(1-s)(2a -1)}\sum_{n=1}^\infty
a_nn^{(1-s)(2a -1)} \int_{2\pi n/N\Delta}^\infty K_{iR}(y)y^{a-s}dy\cr
&=&(2\pi)^{-s}\sum_{m=1}^\infty \left(\sum_{n=1}^m
a_nn^{-s}\right) \int_{2\pi m\Delta}^{2 \pi (m+1)\Delta}
K_{iR}(y)y^{s-1+a}dy \cr
&&\pm (-1)^{a+1} N^{\frac12-s}(2\pi)^{s-1}\sum_{m=1}^\infty
\left(\sum_{n=1}^m a_n n^{s-1}\right) \int_{2\pi m/ N\Delta}^{2\pi
  (m+1)/ N\Delta} K_{iR}(y)y^{a-s}dy
\end{eqnarray}
We use the above expression, truncating the sum over~$m$,
in our calculations.

To decide where to truncate the sum, we use the estimate
\begin{equation}\label{eqn:kdecay}
K_{iR}(y)\approx \sqrt{\frac{\pi}{2y}} e^{-y}e^{\frac{\pi R}{2} },
\end{equation}
for $y>2R$.  We also must assume a bound on the coefficients~$a_n$.
It is conjectured that $|a_n|\le d(n)$, where $d(n)$ is the number
of divisors of~$n$.  That bound is far from being proven, but since
it holds for all coefficients which have been computed, we assume
it in our estimates.  Finally, to isolate $L(s)$ we must divide the
above expression by $(2 \pi)^{-s} G(s+a)$.  This increases the
required precision because $G(s)$ decays rapidly as $\Im(s)$ increases.
However, this will not be an issue because our concern is with
low-lying zeros.

The Bessel function $K_{iR}(y)$ oscillates rapidly as $y\to 0$,
so it is desirable to have the lower limit of integration
as high as possible.  This is accomplished by choosing
$\Delta=1/\sqrt{N}$.  Since $\Delta$ is a free variable, 
we implemented
an error check that looks for differences in $L(s)$ when different
values of $\Delta$ are used, finding good agreement.

We used the above considerations to produce $L$-functions
values with an error of approximately~$10^{-6}$, for
$t$ up to~$\min(20, R)$.  Since we only seek data to produce statistics
involving neighbor spacing and value distributions, this level of
accuracy is more than adequate.  If one wanted to study
many zeros of a single $L$-function then more
sophisticated methods would be needed.

\end{document}